

\documentstyle{amsppt}
\magnification=\magstep1
\NoRunningHeads

\vsize=7.4in


\def\bull{\vrule height .9ex width .8ex depth -.1ex}


\topmatter
\title
The existence of primitives for continuous functions in a quasi-Banach
space
\endtitle

\author
N.J. Kalton \\
University of Missouri-Columbia
\endauthor

\address
Department of Mathematics,
University of Missouri
Columbia, MO  65211, U.S.A.
\endaddress
\email mathnjk\@mizzou1.bitnet \endemail
\thanks Supported by NSF grant DMS-9201357
 \endthanks
\abstract
We show that if $X$ is a quasi-Banach space with
trivial dual then every continuous function $f:[0,1]\to X$ has a
primitive, answering a question of M.M. Popov.
\endabstract

\endtopmatter

\document
\baselineskip=14pt

Let $X$ be a quasi-Banach space and let $f:[0,1]\to X$ be a continuous
function.  We say that $f$ has a primitive if there is a differentiable
function $F:[0,1]\to X$ so that $F'(t)=f(t)$ for $0\le t\le 1.$
M.M. Popov has asked where every continuous function $f:[0,1]\to L_p$
where $0<p<1$ has a primitive; more generally, he asks the same question
for any space with trivial dual \cite4.
 We show here that the answer to this
question is positive.
  We remark that
by an old result of Mazur and Orlicz \cite 3, \cite 6,
every continuous $f$ is Riemann-integrable if and only if $X$ is a
Banach space.

Let us suppose for convenience that $X$ is $p$-normed where $0<p<1$, and
let $I=[0,1].$  Let
$C(I;X)$ be the usual quasi-Banach space of continuous functions
$f:I\to X$ with the quasi-norm $\|f\|_{\infty}=\max_{0\le t\le
1}\|f(t)\|.$  We also introduce the space $C^1(I:X)$ of all functions
$f\in C(I;X)$ which are differentiable at each $t$ and such that the
function
$g:I^2\to X$ is continuous where $g(t,t)=f'(t)$ for $0\le t\le 1$ and
$$ g(s,t)=\frac{f(s)-f(t)}{s-t} $$
when $s\neq t$.  It is easily verified that $C^1(I;X)$ is a quasi-Banach
space under the quasi-norm $$\|f\|_{C^1}=\|f(0)\| + \sup_{0\le s<t\le 1}
\frac{\|f(t)-f(s)\|}{t-s}.$$
Let $C^1_0(I;X)$ be the closed subspace of $C^1(I;X)$ of all $f$ such
that $f(0)=0.$
We consider the map $D:C^1_0(I;X)\to C(I;X)$ given by $Df=f'.$  The
following result is proved in \cite 1.

\proclaim{Theorem 1}If $X$ has trivial dual then for every $x\in X$ there
exists $f\in C^1_0(I;X)$ such that $Df=0$ and $f(1)=x.$\endproclaim

>From this we deduce  the answer to the question of Popov.

\proclaim{Theorem 2}If $X$ has trivial dual then the map
$D:C^1_0(I;X)\to C(I;X)$ is surjective.  In particular every
continuous
$f:I\to X$ has a primitive.\endproclaim

\demo{Proof}From Theorem 1 and the Open Mapping Theorem we deduce the
existence of a constant $M\ge 1$ so that if $x\in X$ there exists $f\in
C^1_0(I;X)$ so that $Df=0,$  $f(1)=x$ and $\|f\|_{C^1}\le M\|x\|.$

Now suppose $g\in C(I;X)$ with $\|g\|_{\infty}< 1.$  For any
$\epsilon>0$ we show the existence of $f\in C^1_0(I;X)$ with
$\|Df-g\|_{\infty}<\epsilon$ and $\|f\|_{C^1}< 4^{1/p}M.$  Once this is
achieved the Theorem follows again from a well-known variant of the Open
Mapping Theorem.

Since $g$ is uniformly continuous, there is a piecewise linear
function $h$ so that $\|g-h\|_{\infty}<\epsilon$ and $\|h\|_{\infty}<1.$
Since $h$ has finite-dimensional range there exists $H\in C^1_0(I;X)$
with $DH=h.$   Now let $n$ be a natural number, and let
$x_{kn}=H(k/n)-H((k-1)/n).$ For
$k=1,2,\ldots n$
define $f_{k,n}\in C^1_0(I;X)$ so that $Df=0,$ $\|f_{k,n}\|_{C^1_0}\le
M\|x_{kn}\|$ and $f_{k,n}(1)=x_{kn}.$ Then we define $F_n\in C^1_0(I;X)$
by
$$ F_n(t) = H(t)-H(\frac{k-1}n)-f_{kn}(nt-k+1)$$
for $(k-1)/n\le t\le k/n.$  Clearly $DF_n=DH=h.$  It remains to estimate
$\|F_n\|_{C^1_0}.$

Let $$\eta(\epsilon)=\sup_{|t-s|\le
\epsilon}\frac{\|H(t)-H(s)\|}{|t-s|}.$$  It is easy to see that
$\lim_{\epsilon\to 0}\eta(\epsilon)=\|h\|_{\infty}<1.$
Now suppose $\frac{k-1}n\le s<t\le \frac kn$ for some $1\le k\le n.$
Then
$$
\align
 \|F_n(t)-F_n(s)\| &\le (\eta(\frac 1n)^p +
n^p\|f_{kn}\|_{C^1}^p)^{1/p}(t-s)
 \\
&\le (\eta(\frac1n)^p + M^pn^p\|x_{kn}\|^p)^{1/p}(t-s)\\
&\le (M^p+1)^{1/p}\eta(\frac 1n)(t-s).
\endalign
$$
Since $F_n(\frac kn)=0$ for $0\le k\le n$ we obtain that for any
$0\le s< t\le 1,$
$$ \|F_n(t)-F_n(s)\| \le 2^{1/p}(M^p+1)^{1/p}\eta(\frac 1n)\min(t-s,\frac
1n).$$

By taking $n$ large enough we have $\|F_n\|_{C^1_0} < 4^{1/p}M.$  Thus
the theorem follows.
\bull\enddemo

We close with a few remarks on the general problem of classifying those
quasi-Banach spaces $X$ so that the map $D:C_0^1(I;X)\to C(I;X)$ is
surjective; let us say that such a space is a $D-$space.
The following facts are clear:

\proclaim{Proposition 3}(1) Any quotient of a D-space is a D-space.
\newline
(2) If $X$ and $Y$ are D-spaces then $X\oplus Y$ is a D-space.
\endproclaim

\demo{Proof}(1) Let $E$ be a closed subspace of $X$ and let $\pi:X\to
X/E$ be the quotient map. Let
$\tilde
\pi:C(I:X)\to C(I;X/E)$ be the induced map $\tilde\pi f= f\circ\pi.$
We start with the observation that $\tilde\pi$ is surjective.
If $g\in C(I;X/E)$ with $\|g\|_{\infty}<1$ then we can find
$f\in
C(I;X)$ with $\|f\|_{\infty}<2^{1/p-1}$ and $\|\tilde\pi f
-g\|_{\infty}<1.$ To do thi suppose $N$ is an integer and let $f_N$ be a
function which is linear on each interval $[(k-1)/N,k/N]$ for $1\le k\le
N$ and such that $\pi f_N(k/N)=g(k/N)$ with $\|f_N(k/N)\|<1$  for $0\le
k\le
N.$  For large enough $N$ we have $\|g-\tilde\pi f_N\|_{\infty}<1$
and our claim is substantiated.

Now if $X$ is a D-space and $g\in C(I;X/E)$ then there exists $f\in
C(I;X)$ with $\tilde\pi f=g.$  Let $F\in C^1_0(I;X)$ with $DF=f.$  Then
if $G=\tilde\pi F$ we have $DG=g.$

(2) is trivial. \bull\enddemo

  In \cite 1 the
notion of the core is defined:  if $X$ is a quasi-Banach space then core
$X$ is the maximal subspace with trivial dual.

\proclaim{Theorem 4} If core $X=\{0\}$ then $X$ is a D-space if and only
if $X$ is a Banach space (i.e. is locally
convex).\endproclaim

\demo{Proof}Suppose core $X=\{0\}$ and $X$ is a D-space.  Suppose
$DF=0$ where $F\in C^1_0(I;X).$  Let $Y$ be the closed subspace generated
by $\{F(s):0\le s\le 1\}.$  We show $Y=\{0\};$  if not there exists a
nontrivial
continuous linear functional $y^*$ on $Y.$  Then $D(y^*\circ F)=0$
so that $y^*(F(s))=0$ for $0\le s\le 1.$  But then $y^*=0$ on $Y.$
We conclude that $Y=\{0\}$ and so $F=0.$  Hence $D$ is one-one and
surjective and by the Closed Graph Theorem $D$ is an isomorphism.

Let $M$ be a constant so that $\|DF\|_{\infty}\le 1$ implies
$\|F\|_{C^1}\le M$ for $F\in C^1_0(I;X).$  Let $\phi$ be any
$C^{\infty}-$ real function on $\bold R$ with $\phi(t)=0$ for $t\le
0$ and $\phi(t)=1$ for $t\ge 1.$
  Let $K=\max_{0\le t\le 1}|\phi'(t)|.$
For any $N$ and any $x_1,\ldots x_N\in X$ with $\max \|x_k\|\le 1,$ we
define
$F(t) =\sum_{k=1}^N \phi(Nt-k+1)x_k$.  Then $F\in C^1_0(I;X)$ and
$\|DF\|_{\infty}\le NK.$  Hence $\|F(1)\| \le NMK,$ i.e.
$$ \|\frac1N(x_1+\cdots+x_N)\| \le MK.$$
This implies $X$ is locally convex.\bull\enddemo

Combining Proposition 3 and Theorem 4 gives that if $X$ is a D-space then
$X/\text{core }X$ is a Banach space.
It is, however, possible to construct an example to show that the
converse to this statement is false, and there does not seem, therefore
to be any nice classification of D-spaces in general.

To construct the example we observe the following theorem.  First for any
quasi-Banach space $X$ let $a_N(X)=\sup\{\|x_1+\cdots+x_N\|:\|x_i\|\le
1\}$ (so that $a_N(X)\ge N).$

\proclaim{Theorem 5}Suppose $X$ is a D-space; then for some constant
$C$ we have $a_N(X)\le Ca_N(\text{core }X).$\endproclaim

\demo{Proof}Let $b_N=a_N(\text{core X}).$   Suppose $x_1,\ldots,x_N\in X$
with $\|x_i\|\le 1$ and define as in Theorem 4,
$F(t)=\sum_{k=1}^N\phi(Nt-k+1)x_k$.  Then $\|DF\|_{\infty}\le NK$ and so
by the Open Mapping Theorem, for some constant $M=M(X)$, there exists
$G\in C_0^1(I;X)$ with $DG=DF$ and $\|G\|_{C^1}\le MNK.$  Then
$\|G(k/N)-G((k-1)/N)\| \le MK$ for $1\le k\le N$.

  Let $H(t)=F(t)-G(t).$  Since $DH=0$ and $X/\text{core X}$ is a Banach
space $H$ has range in core $X$.  Now for $1\le
k\le N,$ $H(k/N)-H((k-1)/N)=x_k
-(G(k/N)-G((k-1)/N)$ so that $\|H(k/N)-H((k-1)/N)\|\le (M^pK^p+1)^{1/p}.$
Hence if $C^p=M^pK^p+1$, we have
$\|H(1)\|\le Cb_N$ or $\|x_1+\cdots+x_n\|\le Cb_n.$\enddemo

To construct our example we start with the Ribe space $Z$
(\cite{2},\cite{5})
which is a space with a one-dimensional subspace $L$ so that $Z/L$ is
isomorphic to $\ell_1$.  A routine calculation shows $a_N(Z)\ge
cN\log N$ for some
$c>0.$ Then let $Y$ be any quasi-Banach space with trivial dual so that
$a_N(Y)=o(N\log N)$  (for example a Lorentz space $L(1,p)$ where
$1<p<\infty$).  Let $j:L\to Y$ be an isometry and let $X$ be the quotient
of $Y\times Z$ by the subspace of all $(jz,z)$ for $z\in L$.  Then $Z$
embeds into $X$ so that $a_n(X)\ge cN\log N$ but $\text{core X}\sim Y$ so
that $X$ cannot be a D-space.  However $X/\text{core X}$ is isomorphic to
$Z/L$ which is a Banach space.

\enddocument